\documentclass[12pt,twoside,a4paper]{article}
\usepackage{amsmath,amssymb}
\usepackage{fullpage}

\usepackage{pslatex}
\usepackage{url}
\usepackage{comment}

\newcommand{\Odip}[2]{\mathcal{O}_{#1}\!\left(#2\right)\mathchoice{\!}{}{}{}}
\newcommand{\Odi}[1]{\Odip{}{#1}}
\newcommand{\Eulerphi}{\varphi}
\newcommand{\dx}{\mathrm{d}}
\newcommand{\Zs}[1]{\mathbb{Z}_{#1}^*}
\newcommand{\lcm}{\operatorname{lcm}}
\newcommand{\A}{\mathcal{A}}
\newcommand{\B}{\mathcal{B}}
\newcommand{\Q}{\mathbb{Q}}
\newcommand{\Res}{\operatornamewithlimits{Res}}

\newtheorem{Theorem}{Theorem}
\newtheorem{Lemma}{Lemma}
\newenvironment{Proof}[1][Proof]{\par\noindent\textbf{#1.}~}
  {\hfill$\square$\smallskip\par}

\begin{document}

\title{On the constant in the Mertens product \\
       for arithmetic progressions. \\
       I. Identities}
\author{A.~LANGUASCO and A.~ZACCAGNINI}

\date{September 26, 2008}

\maketitle

\begin{abstract}
We prove new identities for the constant in the Mertens product over
primes in the arithmetic progressions $a \bmod q$.
\end{abstract}

\section{Introduction}

Let $a$, $q$ be integers with $(q, a) = 1$ and denote by $p$ a
prime number.
In 1974 Williams \cite{Williams74} proved that
\begin{equation}
\label{def-C}
  P(x; q, a)
  =
  \prod_{\substack{p \le x \\ p \equiv a \bmod q}}
    \Bigl( 1 - \frac1p \Bigr)
   =
  \frac{C(q, a)}{(\log x)^{1/\Eulerphi(q)}}
  +
  \Odi{ \frac1{(\log x)^{1/\Eulerphi(q)+1} } }
\end{equation}
as $x \to +\infty$, where $C(q, a)$ is real and positive and satisfies
\[
  C(q,a)^{\Eulerphi(q)}
  =
  e^{-\gamma} \frac{q}{\Eulerphi(q)}
  \prod_{\chi\neq \chi_0}
  \Big(
    \frac{K(1,\chi)}{L(1,\chi)}
  \Big)^{\overline{\chi}(a)},
\]
where $\gamma$ is the Euler constant, $\Eulerphi$ is the Euler totient
function, $L(s, \chi)$ is the Dirichlet $L$-function associated to the
Dirichlet character $\chi \bmod q$ and $\chi_0$ is the principal
character to the modulus $q$.
The function $K$ is defined by means of
\[
  K(s, \chi)
  =
  \sum_{n=1}^{+\infty} k_{\chi}(n) n^{-s},
\]
where $k_{\chi}(n)$ is the completely multiplicative function whose
value at primes is given by
\[
  k_{\chi}(p)
  =
  p \left( 1 - \Big( 1 - \frac{\chi(p)}{p} \Big)
               \Big( 1 - \frac{1}{p} \Big)^{-\chi(p)}
    \right).
\]
In our recent paper \cite{LanguascoZaccagnini2007} we obtained a
version of Williams's result stated in \eqref{def-C} which is uniform
in the $q$ aspect.
In the same paper, as a by-product, we also obtained the following
elementary expression for $C(q,a)$:
\begin{equation}
\label{value-C}
  C(q, a)^{\Eulerphi(q)}
  =
  e^{-\gamma}
  \prod_p
    \Bigl( 1 - \frac1p \Bigr)^{\alpha(p; q, a)}
\end{equation}
where $\alpha(p; q, a) = \Eulerphi(q) - 1$ if $p \equiv a \bmod q$ and
$\alpha(p; q, a) = -1$ otherwise.
The infinite product is convergent, though not absolutely, by the
Prime Number Theorem for Arithmetic Progressions.
Actually, a slightly simpler proof of identity \eqref{value-C} than
the one we gave in \cite{LanguascoZaccagnini2007}, sect.~6, can be
obtained as follows, once the relevant limit is known to exist: taking
logarithms in \eqref{def-C} and using the classical Mertens Theorem,
we find
\begin{align*}
  \log C(q, a)
  &=
  \lim_{x \to +\infty}
    \Bigl\{
      \log
      \prod_{\substack{p \le x \\ p \equiv a \bmod q}}
        \Bigl( 1 - \frac1p \Bigr)
      +
      \frac1{\Eulerphi(q)}
      \log \log x
    \Bigr\} \\
  &=
  -\frac{\gamma}{\Eulerphi(q)}
  +
  \lim_{x \to +\infty}
    \Bigl\{
      \log
      \prod_{\substack{p \le x \\ p \equiv a \bmod q}}
        \Bigl( 1 - \frac1p \Bigr)
      -
      \log \prod_{p \le x}
        \Bigl( 1 - \frac1p \Bigr)^{1 / \Eulerphi(q)}
    \Bigr\} \\
  &=
  -\frac{\gamma}{\Eulerphi(q)}
  +
  \frac1{\Eulerphi(q)}
  \log
  \Bigl\{
    \lim_{x \to +\infty}
      \prod_{p \le x}
        \Bigl( 1 - \frac1p \Bigr)^{\alpha(p; q, a)}
  \Bigr\}.
\end{align*}
The product in \eqref{value-C} is very slowly convergent and it is
difficult to compute an accurate numerical approximation to $C(q, a)$
from it.
Our aim here to give a different form for the constant defined in
\eqref{value-C}: unfortunately, this form is not suitable for
numerical computations, a problem we tackle in part II
\cite{LanguascoZaccagnini2008a}.

As a corollary of the identities proved in the first part of the
paper, we derive the formulae that Uchiyama \cite{Uchiyama71}
gave in the case $q = 4$ and $a \in \{1$, $3\}$, though Uchiyama's
direct proof is obviously much simpler.
We also derive the explicit expressions that Williams gave in
Theorem~2 of \cite{Williams74} for $C(24, a)$ for every integer $a$
such that $(24, a) = 1$ and the ones that Grosswald \cite{Grosswald87}
obtained for $C(q, a)$ for $q \in \{4$, $6$, $8\}$ and every integer
$a$ coprime to $q$. 
We recall that, in Proposition~1 of \cite{Moree2006}, Moree gives formula
\eqref{value-c-cyclic} below when $q$ is a prime number and $a = 1$.

The statement of these new formulae themselves is not simple, and
reflects both the structure of the group $\Zs{q}$ and the properties
of the residue class $a$.
For this reason, we will not state a formal theorem here, but rather
point to the various results as formulae \eqref{value-c-cyclic} or its
alternative version \eqref{value-c-cyclic-alt} when $\Zs{q}$ is cyclic
and $a = 1$, then \eqref{value-c-general-alt} for general $q$ and $a = 1$,
and \eqref{value-c-more-general} in the most general case.
We think our results will be clearer if we go through stages of
increasing generality.

We may summarize our main result saying that, in the case $a = 1$, for
each reduced residue class $b$ we have to determine a positive integer
$t_b$ (actually, the order of $b$ in the multiplicative group
$\Zs{q}$) and we will then express $C(q, 1)$ as a sort of Euler
product where a prime $p$ has the exponent $-t_p$.
Collecting all residue classes of maximal order, we may reduce the
number of factors needed, at the price of the computation of the power
of a suitable value of the Riemann zeta function at an even integer.
The important feature of \eqref{value-c-general-alt} is that $t_b \ge 2$
for all $b \ne 1$, whereas the exponent of all primes in the Euler
factors in \eqref{value-C} are $-1$.

The case $a \ne 1$ is genuinely more complicated: in fact, in general
it will not be possible to give a simple closed form for the Euler
factors, though they can be expressed by means of a rapidly convergent
power series.
On the other hand, the constant defined in \eqref{Meissel} will arise:
it is related to the so-called Meissel--Mertens constant; for its
computation in the case $q = 1$ see \S2.2 of Finch \cite{Finch2003}.

We would like to thank Pieter Moree for providing us some references
and Giuseppe Molteni for suggesting a simplification in the proof of
Lemma~\ref{Lemma-1}.

\section{Reduction to character sums}

It turns out to be better to get rid of the prime factors of $q$ at
the outset: therefore, we let $c(q, a)$ be defined by means of
\begin{equation}
\label{def-c}
  C(q, a)^{\Eulerphi(q)}
  =
  e^{-\gamma}
  \frac q{\Eulerphi(q)}
  \,
  c(q, a)
  \qquad\text{so that}\qquad
  c(q, a)
  =
  \prod_{p \nmid q}
    \Bigl( 1 - \frac1p \Bigr)^{\alpha(p; q, a)}
\end{equation}
and let
\[
  c(x; q, a)
  =
  \prod_{\substack{p \le x \\ p \nmid q}}
    \Bigl( 1 - \frac1p \Bigr)^{\alpha(p; q, a)}
\]
denote its partial product.
Our strategy is to express the product of $c(x; q, a)$ and partial
products of powers of $L(1, \chi)$ as a complicated but quickly
convergent product, where $\chi$ ranges over all non-principal
Dirichlet characters modulo $q$.
When necessary, we use the abbreviation
\begin{equation}
\label{prod-L(1)}
  \Pi(q, a)
  =
  \prod_{\substack{\chi \bmod q \\ \chi \ne \chi_0}}
    L(1, \chi)^{-\overline{\chi}(a)}.
\end{equation}
By orthogonality, we have
\begin{equation}
\label{main-prod}
  c(x; q, a)
  \prod_{\substack{\chi \bmod q \\ \chi \ne \chi_0}}
    \,\,
    \prod_{p \le x}
      \Bigl( 1 - \frac{\chi(p)}p \Bigr)^{-\overline{\chi}(a)}
  =
  \prod_{\substack{p \le x \\ p \nmid q}}
  \,\,
  \prod_{\substack{\chi \bmod q \\ \chi \ne \chi_0}}
    \Bigl\{
      \Bigl( 1 - \frac1p         \Bigr)^{ \overline{\chi}(a)\chi(p)}
      \Bigl( 1 - \frac{\chi(p)}p \Bigr)^{-\overline{\chi}(a)}
    \Bigr\}.
\end{equation}
Using the Taylor series expansion of $\log(1 - t)$ we see that
\begin{equation}
\label{def-S}
  \sum_{\substack{\chi \bmod q \\ \chi \ne \chi_0}}
  \log
  \Bigl\{
    \Bigl( 1 - \frac1p         \Bigr)^{ \overline{\chi}(a)\chi(p)}
    \Bigl( 1 - \frac{\chi(p)}p \Bigr)^{-\overline{\chi}(a)}
  \Bigr\}
  =
  \sum_{m \ge 2} \frac1{m p^m} \, S_m(p; q, a),
\end{equation}
say, where $S_m(p; q, a)$ is the character sum defined by
\begin{equation}
\label{def-Sm}
  S_m(p; q, a)
  =
  \sum_{\chi \bmod q}
    \overline{\chi}(a) \bigl( \chi^m(p) - \chi(p) \bigr).
\end{equation}
It is this character sum that reflects the structure of $\Zs{q}$ and
the properties of the element $a$.
We will prove below in \eqref{value-S-a=b} and \eqref{value-S-a-ne-b}
that either $p \equiv a \bmod q$, or $S_m(p; q, a)$ vanishes unless
$m$ belongs to a suitable arithmetic progression modulo a divisor of
$\Eulerphi(q)$.
The simplest case, not surprisingly, is when $\Zs{q}$ is cyclic and
$a = 1$.

\section{The character sum $S_m$ in the simplest case}

We notice that, obviously, $S_m(1; q, a) = 0$, and we may assume that
$p \not\equiv 1 \bmod q$.
For the time being, we also assume that $a = 1$.
Let $t_p$ denote the order of $p$ in the multiplicative group
$\Zs{q}$, that is, the smallest positive integer $k$ such that
$p^k \equiv 1 \bmod q$, and notice that $t_p \ge 2$ since
$p \not\equiv 1 \bmod q$.
It is then quite easy to see that
\[
  \sum_{\chi \bmod q}
    \chi^m(p)
  =
  \sum_{\chi \bmod q}
    \chi(p^m)
  =
  \begin{cases}
    \Eulerphi(q) & \text{if $t_p \mid m$} \\
    0            & \text{otherwise}.
  \end{cases}
\]
Hence, using orthogonality and the Taylor series for $\log(1 - t)$
again, we have
\[
  \sum_{m \ge 2} \frac1{m p^m} S_m(p; q, 1)
  =
  \sum_{n \ge 1}
    \frac{\Eulerphi(q)}{n t_p p^{n t_p}}
  =
  \log \Bigl(1 - \frac1{p^{t_p}}\Bigr)^{- \Eulerphi(q) / t_p}.
\]
We classify primes according to their residue class $b \bmod q$, and
notice that $t_p$ depends only on $b$, if $p \equiv b \bmod q$.
Substituting into \eqref{main-prod} and letting $x \to +\infty$, we see
that
\begin{equation}
\label{value-c-cyclic}
  c(q, 1)
  =
  \Pi(q, 1)
  \,\,
  \prod_{b \in \Zs{q} \setminus \{ 1 \}}
  \,\,
  \prod_{p \equiv b \bmod q}
    \Bigl(1 - \frac1{p^{t_b}}\Bigr)^{- \Eulerphi(q) / t_b}.
\end{equation}
We notice that the quantity $\Pi(q, 1)$ is connected to the Dedekind
zeta function of the $q$-th cyclotomic field $K = \Q(\zeta_q)$ by
means of the relation
\[
  \Pi(q, 1)^{-1}
  =
  \Res_{s = 1} \zeta_K(s)
  \prod_{\substack{\chi \bmod q \\ \chi \ne \chi_0}}
    \prod_{p \mid q} \Bigl( 1 - \frac{\chi_f(p)}p \Bigr),
\]
where $\chi_f$ denotes the primitive character that induces $\chi$ and
$f$ is its conductor.

Assume that $\Zs{q}$ is cyclic.
Relation \eqref{value-c-cyclic} is our first formula, and it is worth
noticing that a slightly better form, from the point of view of the
explicit computation of $C(q, 1)$, can be given grouping the
contribution of primes of maximal order.
More specifically, we may rewrite \eqref{value-c-cyclic} as
\allowdisplaybreaks
\begin{align}
  c(q, 1)
  &=
  \Pi(q, 1)
  \prod_{\substack{n \mid \Eulerphi(q) \\ n > 1}} \,\,
  \prod_{\substack{p \phantom{\mid} \\ t_p = n}}
    \Bigl(1 - \frac1{p^n}
    \Bigr)^{-\Eulerphi(q) / n} \notag \\
  &=
  \Pi(q, 1)
  \prod_{\substack{p \phantom{\mid} \\ t_p = \Eulerphi(q)}}
    \Bigl(1 - \frac1{p^{\Eulerphi(q)}} \Bigr)^{-1}
  \prod_{\substack{n \mid \Eulerphi(q) \\ 1 < n < \Eulerphi(q)}} \,\,
  \prod_{\substack{p \phantom{\mid} \\ t_p = n}}
    \Bigl(1 - \frac1{p^n}
    \Bigr)^{-\Eulerphi(q) / n} \notag \\
  &=
  \Pi(q, 1)
  \prod_{p \nmid q}
    \Bigl(1 - \frac1{p^{\Eulerphi(q)}} \Bigr)^{-1}
  \prod_{p \equiv 1 \bmod q}
    \Bigl(1 - \frac1{p^{\Eulerphi(q)}} \Bigr)
  \notag \\
  &\qquad\qquad\times
  \prod_{\substack{n \mid \Eulerphi(q) \\ 1 < n < \Eulerphi(q)}} \,\,
  \prod_{\substack{p \phantom{\mid} \\ t_p = n}}
    \Bigl\{
      \Bigl(1 - \frac1{p^n}
      \Bigr)^{-\Eulerphi(q) / n}
      \Bigl(1 - \frac1{p^{\Eulerphi(q)}} \Bigr)
    \Bigr\} \notag \\
  &=
  \zeta(\Eulerphi(q)) \,
  \Pi(q, 1)
  \prod_{p \mid q} \Bigl(1 - \frac1{p^{\Eulerphi(q)}} \Bigr)
  \prod_{p \equiv 1 \bmod q}
    \Bigl(1 - \frac1{p^{\Eulerphi(q)}} \Bigr)
  \notag \\
\label{value-c-cyclic-alt}
  &\qquad\qquad\times
  \prod_{\substack{n \mid \Eulerphi(q) \\ 1 < n < \Eulerphi(q)}} \,\,
  \prod_{\substack{p \phantom{\mid} \\ t_p = n}}
    \Bigl\{
      \Bigl(1 - \frac1{p^n}
      \Bigr)^{-\Eulerphi(q) / n}
      \Bigl(1 - \frac1{p^{\Eulerphi(q)}} \Bigr)
    \Bigr\}.
\end{align}
The value of $\zeta(\Eulerphi(q))$ is easily computed, at least when
$q$ is comparatively small, by means of the Bernoulli numbers, since
$\Eulerphi(q)$ is even for $q \ge 3$.
It is also worth noticing that the exponents of the prime numbers in
the last product are all at least 2, though they are usually much
larger.
Uchiyama's formula \eqref{Uchiyama} for $C(4, 1)$ in \cite{Uchiyama71}
is the case $q = 4$ of the above expression: there is only one
non-principal character $\chi$ modulo $4$, and $L(1, \chi) = \pi / 4$.
The last product is empty, and $\zeta(2) = \pi^2 / 6$.
The formula for $C(4, 3)$ is easily deduced from this and the
classical Mertens Theorem since $C(4, 3) C(4, 1) = 2 e^{-\gamma}$.
Grosswald's formula \eqref{Grosswald1} for $C(6, 1)$ in
\cite{Grosswald87} is another special case of \eqref{def-c} and
\eqref{value-c-cyclic-alt} since $\Eulerphi(6) = 2$, there is only one
non-principal character $\chi$ modulo~$6$, and
$L(1, \chi) = \pi / (2 \sqrt{3})$.
By the Mertens Theorem we have $C(6, 5) C(6, 1) = 3 e^{-\gamma}$ and
the formula for $C(6, 5)$ can be easily deduced.
Our formula \eqref{value-c-cyclic-alt} also contains Moree's
\cite{Moree2006} which is the special case where $q$ is prime.

\section{The sum $S_m$ in the general case}

\begin{Lemma}
\label{Lemma-1}
Let $S_m$ be the character sum defined in \textup{\eqref{def-Sm}}.
If $a \equiv b \bmod q$ then there exists a positive integer $t_a$
dividing $\Eulerphi(q)$ such that
\begin{equation}
\label{value-S-a=b}
  S_m(a; q, a)
  =
  \begin{cases}
    -\Eulerphi(q)
      & \text{if $m \not\equiv 1 \bmod t_a$,} \\
    0
      & \text{otherwise.}
  \end{cases}
\end{equation}
In this case, $t_a$ is precisely the order of $a$ in the group
$\Zs{q}$, so that $t_a \ge 2$ unless $a = b = 1$.
If $a \not\equiv b \bmod q$ then either the equation
$b^y \equiv a \bmod q$ has no solution $y \in \mathbb{N}$, or there
exist a positive integer $t_b$ dividing $\Eulerphi(q)$ and an integer
$s_b$ such that $1 \le s_b \le t_b$ and
\begin{equation}
\label{value-S-a-ne-b}
  S_m(b; q, a)
  =
  \begin{cases}
    \Eulerphi(q)
      & \text{if $b^y \equiv a \bmod q$ has a solution, and
              $m \equiv s_b \bmod t_b$,} \\
    0
      & \text{otherwise.}
  \end{cases}
\end{equation}
In particular, if $a = 1$ then $S_m(1; q, 1) = 0$, while, if
$b \not\equiv 1 \bmod q$, then there is an integer $t_b \ge 2$ such
that $S_m(b; q, 1) = \Eulerphi(q)$ if $m \equiv 0 \bmod t_b$ and is
$0$ otherwise.
\end{Lemma}

\begin{Proof}
The complete multiplicativity of the Dirichlet characters implies that
\[
  S_m(b; q, a)
  =
  \sum_{\chi \bmod q}
    \bigl( \chi(b^m a^{-1}) - \chi(b a^{-1}) \bigr)
\]
If $a \equiv b \equiv 1 \bmod q$ then $S_m(a; q, a) = 0$.
If $a \equiv b \not\equiv 1 \bmod q$ then the equation
$a^y \equiv a \bmod q$ has the solution $y \equiv 1 \bmod t_a$, where
$t_a$ is the order of $a$ in $\Zs{q}$.
Hence, $S_m(a; q, a) = 0$ if $m \equiv 1 \bmod t_a$ and
$S_m(a; q, a) = -\Eulerphi(q)$ otherwise, by orthogonality.

If $a \not\equiv b \bmod q$ and the equation $b^y \equiv a \bmod q$
has no solution, then $S_m(b, q, a) = 0$ by orthogonality.
If the equation above has the solution $y \equiv s_b \bmod t_b$ (where
$t_b \ge 2$ is a suitable divisor of $\Eulerphi(q)$ and $s_b$ is an
integer with $1 \le s_b \le t_b$) then $S_m(b; q, a) = \Eulerphi(q)$
if $m \equiv s_b \bmod t_b$ and is $0$ otherwise.
\end{Proof}

It is quite clear from the Lemma above that the case $a = 1$ is
indeed much simpler than the general one.
We see that we have proved that \eqref{value-c-cyclic} holds also in
the general case, where $t_b$ is the divisor of $\Eulerphi(q)$ such
that the solution of the equation $b^y \equiv 1 \bmod q$ is the class
$0 \bmod t_b$, that is, the order of $b$ in $\Zs{q}$.
The computations that lead to \eqref{value-c-cyclic-alt} are still
essentially valid, with one modification.
We recall the definition of the Carmichael $\lambda$ function: if
$q = 2^\alpha p_1^{\alpha_1} \cdots p_k^{\alpha_k}$, then
\begin{equation}
\label{def-lambda}
  \lambda(q)
  =
  \begin{cases}
    \lcm\{
      \Eulerphi(2^\alpha), \Eulerphi(p_1^{\alpha_1}), \dots,
      \Eulerphi(p_k^{\alpha_k})
    \}
    & \text{if $\alpha \le 2$;} \\
    \lcm\{
      2^{\alpha - 2}, \Eulerphi(p_1^{\alpha_1}), \dots,
      \Eulerphi(p_k^{\alpha_k})
    \}
    & \text{if $\alpha \ge 3$.}
  \end{cases}
\end{equation}
In other words, $\lambda(q)$ is the highest order of the elements of
the group $\Zs{q}$.
Let
$\A(q) = \{ b \in \Zs{q} \setminus \{ 1 \} \colon b^k \equiv 1 \bmod q$
for some positive $k < \lambda(q) \}$ denote the set of elements of
$\Zs{q} \setminus \{ 1 \}$ whose order is not maximal.
Arguing as in the proof of \eqref{value-c-cyclic-alt}, that is,
grouping the contribution of the primes of maximal order, we obtain
the following Theorem, which generalizes Williams's \cite{Williams74}
and Grosswald's \cite{Grosswald87} formulae for $C(24, 1)$ and $C(8, 1)$
respectively and contains \eqref{value-c-cyclic-alt} as a special case.

\begin{Theorem}
For all integers $q \ge 3$ the value of the constant $c(q, 1)$ is
given by
\begin{align}
  c(q, 1)
  &=
  \zeta(\lambda(q))^{\Eulerphi(q)/\lambda(q)} \,\,
  \Pi(q, 1)
  \,\,
  \prod_{p \mid q}
    \Bigl( 1 - \frac1{p^{\lambda(q)}} \Bigr)^{\Eulerphi(q)/\lambda(q)}
  \prod_{p \equiv 1 \bmod q}
    \Bigl( 1 - \frac1{p^{\lambda(q)}} \Bigr)^{\Eulerphi(q)/\lambda(q)}
  \notag \\
\label{value-c-general-alt}
  &\qquad\qquad \times
  \prod_{b \in \A(q)}
  \,\,
  \prod_{p \equiv b \bmod q}
    \Bigl\{
      \Bigl(1 - \frac1{p^{t_b}}\Bigr)^{-\Eulerphi(q) / t_b}
      \Bigl( 1 - \frac1{p^{\lambda(q)}} \Bigr)^{\Eulerphi(q)/\lambda(q)}
    \Bigr\}
\end{align}
where $\Pi(q, a)$ is defined in \textup{\eqref{prod-L(1)}}, $\lambda$ is
the Carmichael lambda function defined in \textup{\eqref{def-lambda}},
$t_b$ denotes the order of $b$ in the group $\Zs{q}$ and
$\A(q) = \{ b \in \Zs{q} \setminus \{ 1 \} \colon t_b < \lambda(q) \}$.
\end{Theorem}

\section{The general formula}

In this section we assume that $a \not\equiv 1 \bmod q$, and we let
$\B(q, a)$ denote the set $\{ b \in \Zs{q} \setminus \{1, a\} \colon$
the equation $b^y \equiv a \bmod q$ has a solution$\}$.
For the elements of this set, we implicitly define the integers $t_b$
and $s_b$ as in the proof of Lemma~\ref{Lemma-1}.
In order to state the formula corresponding to
\eqref{value-c-general-alt} in general, we need to introduce the
function $f_{t,s}$ which is defined, for positive integers $t$ and $s$
with $1 \le s \le t$ and real $x$ with $|x| < 1$, by means of the
relation
\begin{equation}
\label{def-f}
  f_{t,s}(x)
  =
  \sum_{\substack{n \ge 1 \\ n \equiv s \bmod t}}
    \frac{x^n}n
  =
  \int_0^x \frac{u^{s - 1}}{1 - u^t} \, \dx u.
\end{equation}
The rightmost equality is proved computing the derivative of the
function $f_{t,s}$ and then summing the ensuing geometric progression.
Notice that, when $s = t$, a closed form for the integral can be
easily given in terms of the logarithmic function, as we did above:
indeed, $f_{t,t}(x) = - t^{-1} \log(1 - x^t)$.
In general, the Taylor series for $f_{t,s}$ is fairly quickly
convergent since we will compute it at $x = p^{-1}$.

Lemma~\ref{Lemma-1} above amounts to saying that, given integers $q$
and $a$ and a reduced residue class $b \bmod \Eulerphi(q)$, either the
equation $b^y \equiv a \bmod q$ does not have a solution, or its
solution is a congruence class $s_b \bmod t_b$, where
$t_b \mid \Eulerphi(q)$ and we may assume that $1 \le s_b \le t_b$ and
that $t_b \ge 2$, since $t_b = 1$ if and only if $a = b = 1$.
Moreover, unless $a = b$, the congruence class $s_b$ will not be
$1 \bmod q$.
Therefore, classifying primes according to their residue class
$b \bmod q$ again, and substituting either \eqref{value-S-a=b} or
\eqref{value-S-a-ne-b} into \eqref{def-S}, we see that the
corresponding factor in the product \eqref{main-prod} is $1$ if
$b = 1$ or $b \not\in \B(q, a) \cup \{ a \}$, and is
\[
  \prod_{p \equiv b \bmod q}
    \exp\Bigl( \Eulerphi(q) f_{t_b, s_b} (p^{-1}) \Bigr)
\]
if $b \in \B(q, a)$.
Recall that for $b \equiv a \bmod q$ we have $s_a = 1$: hence, for
primes $p \equiv a \bmod q$ we have a factor
\begin{align*}
  \prod_{p \equiv a \bmod q}
    \exp\Bigl(
      -\Eulerphi(q)
      \sum_{\substack{m \ge 2 \\ m \not\equiv 1 \bmod t_a}}
        \frac1{m p^m}
    \Bigr)
  &=
  \prod_{p \equiv a \bmod q}
    \exp\Bigl(
      \Eulerphi(q)
      \sum_{\substack{m \ge 2 \\ m \equiv 1 \bmod t_a}}
        \frac1{m p^m}
      -\Eulerphi(q)
      \sum_{m \ge 2} \frac1{m p^m}
    \Bigr) \\
  &=
  \prod_{p \equiv a \bmod q}
    \exp\Bigl(
      \Eulerphi(q) f_{t_a, 1} \Bigl( \frac1p \Bigr)
      -\Eulerphi(q)
      \sum_{m \ge 1} \frac1{m p^m}
    \Bigr) \\
  &=
  \exp\Bigl(
    \Eulerphi(q)
    \sum_{p \equiv a \bmod q}
      \Bigl( \frac1p + \log \Bigl(1 - \frac1p \Bigr) \Bigr)
    \Bigr) \\
  &\qquad\qquad\times
  \exp\Bigl(
    \Eulerphi(q)
    \sum_{p \equiv a \bmod q}
      \Bigl( f_{t_a, 1} \Bigl( \frac1p \Bigr) - \frac1p \Bigr)
    \Bigr).
\end{align*}
Finally, collecting all identities, we see that we have proved the
following result.

\begin{Theorem}
For all integers $q \ge 3$ and all integers $a$ such that $(q, a) = 1$
and $a \not\equiv 1 \bmod q$, the value of the constant $c(q, a)$ is
given by
\begin{align}
  c(q, a)
  &=
  \Pi(q, a) \,\,
  \prod_{b \in \B(q, a)}
  \,\,
  \prod_{p \equiv b \bmod q}
    \exp\Bigl( \Eulerphi(q) f_{t_b, s_b} \Bigl( \frac1p \Bigr) \Bigr) \notag \\
\label{value-c-more-general}
  &\qquad\times
  \exp( \Eulerphi(q) B(q, a) )
  \exp\Bigl(
    \Eulerphi(q)
    \sum_{p \equiv a \bmod q}
      \Bigl( f_{t_a, 1} \Bigl( \frac1p \Bigr) - \frac1p \Bigr)
    \Bigr),
\end{align}
where $\B(q, a) = \{ b \in \Zs{q} \setminus \{1, a\} \colon$
the equation $b^y \equiv a \bmod q$ has a solution$\}$,
\begin{equation}
\label{Meissel}
  B(q, a)
  =
  \sum_{p \equiv a \bmod q}
    \Bigl( \frac1p + \log \Bigl(1 - \frac1p \Bigr) \Bigr),
\end{equation}
$\Pi(q, a)$ is defined in \textup{\eqref{prod-L(1)}} and $f_{t, s}$ is
defined in \textup{\eqref{def-f}}.
\end{Theorem}

For the special cases $q \in \{4$, $6$, $8$, $24\}$ and $(q, a) = 1$
with $a \not\equiv 1 \bmod{q}$, equation \eqref{value-c-more-general}
collapses to the formulae given by Uchiyama \cite{Uchiyama71},
Williams \cite{Williams74} and Grosswald \cite{Grosswald87}.

For some special values of $t$ and $s$ it is possible to compute a
closed form for $f_{t, s}$ as in the previous sections, and give a
more explicit result: the following section contains some examples.

\section{Explicit values}

Using \eqref{def-c}, \eqref{value-c-cyclic-alt} and
\eqref{value-c-general-alt}, we can compute explicitly a few values of
the constant $C(q, a)$.
For the evaluation of $L(1, \chi)$ needed to determine $\Pi(q, 1)$ we
refer to Corollary 10.3.2 and Proposition 10.3.5 of Cohen
\cite{Cohen2007b}.
The value \eqref{Uchiyama} is due to Uchiyama while the values
\eqref{Grosswald1} and \eqref{Grosswald2} are due to Grosswald.
Notice that $\A(8) = \varnothing$ and that
$\Pi(8, 1)= 32 \pi^{-2} (\log(3 + 2 \sqrt{2}))^{-1}$.
\begin{align}
\label{Uchiyama}
  C(4, 1)^2
  &=
  \pi \,\, e^{-\gamma}
  \prod_{p \equiv 1 \bmod 4} \Bigl( 1 - \frac1{p^2} \Bigr) \\
\label{Grosswald1}
  C(6, 1)^2
  &=
  \frac{2 \pi \sqrt{3}}3 \,\, e^{-\gamma}
  \prod_{p \equiv 1 \bmod 6} \Bigl( 1 - \frac1{p^2} \Bigr) \\
\label{Grosswald2}
  C(8, 1)^4
  &=
  \frac1{32} \,\, \pi^4 \,\, e^{-\gamma} \,\, \Pi(8, 1)
  \prod_{p \equiv 1 \bmod 8} \Bigl( 1 - \frac1{p^2} \Bigr)^2
  =
  \frac{\pi^2 e^{-\gamma}}{\log(3 + 2\sqrt{2})}
  \prod_{p \equiv 1 \bmod 8} \Bigl( 1 - \frac1{p^2} \Bigr)^2.
\end{align}
For $q = 24$ we recover the value given on page~357 of Williams
\cite{Williams74}: using \eqref{def-c} and \eqref{value-c-general-alt}
we get
\[
  C(24, 1)^8
  =
  \frac{2 \pi^4 e^{-\gamma}}
       {9 \log(2 + \sqrt{3}) \log(1 + \sqrt{2}) \log(5 + 2\sqrt{6})}
  \prod_{p \equiv 1 \bmod 24} \Bigl( 1 - \frac1{p^2} \Bigr)^4
\]
since $\Eulerphi(24) = 8$, $\lambda(24) = 2$, $\A(24) = \varnothing$
and
\[
  \Pi(24, 1)
  =
  \frac{486}
  {\pi^{4} \log(2+\sqrt{3})\log(1+\sqrt{2})\log(5+2\sqrt{6})}.
\]
For the more complicated case $q = 15$ we get
\[
  C(15, 1)^8
  =
  \frac{15}8 \,\,
  \frac{\pi^8}{90^2} \,\, \frac{3328^2}{3375^2} \,\, e^{-\gamma} \,\,
  \Pi(15, 1)
  \prod_{p \equiv 1 \bmod 15} \Bigl( 1 - \frac1{p^4} \Bigr)^2
  \prod_{\substack{b \in \{4, 11, 14\} \\ p \equiv b \bmod 15}}
    \Bigl( \frac{1 + p^{-2}}{1 - p^{-2}} \Bigr)^2
\]
where we have used the fact that $\A(15) = \{4$, $11$, $14\}$ and the
values $\zeta(4) = \pi^4 / 90$ and $\lambda(15) = 4$.
A fairly lengthy computation reveals that
\[
  \Pi(15, 1)^{-1}
  =
  \frac{2^9 \, \pi^4}{3 \cdot 15^4}
  \Bigl( \log^2\Bigl( \frac{\sin(  \pi / 15)}{\sin(4 \pi / 15)} \Bigr)
         +
	 \log^2\Bigl( \frac{\sin(7 \pi / 15)}{\sin(2 \pi / 15)} \Bigr)
  \Bigr)
  \log\Bigl( \frac{1 + \sqrt{5}}2 \Bigr).
\]
This can be made more explicit using suitable trigonometrical
identities and the value for $\sin(\pi / 15)$,

Using \eqref{value-c-cyclic-alt} we can compute $C(5, 1)$:
\begin{align*}
  C(5, 1)^4
  &=
  \frac{5}{4}
  \,\frac{\pi^4}{90} \,\, e^{-\gamma} \,\,
  \Pi(5, 1)
  \,\frac{624}{625}
  \prod_{p \equiv 1 \bmod 5} \Bigl( 1 - \frac1{p^4} \Bigr)
  \prod_{p \equiv 4 \bmod 5} \Bigl( 1 - \frac1{p^2} \Bigr)^{-2}
  			\Bigl( 1 - \frac1{p^4} \Bigr) \\
  &=
  \frac{13 \sqrt{5}\,\, \pi^2\,\, e^{-\gamma} }{150\,\, \log((1+\sqrt{5})/2)}
  \prod_{p \equiv 1 \bmod 5} \Bigl( 1 - \frac1{p^4} \Bigr)
  \prod_{p \equiv 4 \bmod 5} \Bigl( \frac{1 + p^{-2}}{1 - p^{-2}} \Bigr)
\end{align*}
where we have used the value
$\Pi(5, 1) = 25 \sqrt{5} \pi^{-2} (4 \log((1+\sqrt{5})/2))^{-1}$ and the
fact that $\A(5) = \{ 4 \}$.

The last examples are for $q = 5$ and $a \in \{2$, $3$, $4\}$.
A short computation shows that
\begin{align*}
  f_{4, 1}(x)
  &=
  \frac14 \log \Bigl( \frac{1 + x}{1 - x} \Bigr)
  +
  \frac12 \arctan(x) \\
  f_{4, 2}(x)
  &=
  \frac14 \log \Bigl( \frac{1 + x^2}{1 - x^2} \Bigr) \\
  f_{4, 3}(x)
  &=
  \frac14 \log \Bigl( \frac{1 + x}{1 - x} \Bigr)
  -
  \frac12 \arctan(x).
\end{align*}
Furthermore, using the fact that $\Zs{5}$ is generated by $2$, we see
that $\B(5, 2) = \{3\}$ (with $s_3 = 3$ and $t_3 = 4$),
$\B(5, 3) = \{2\}$ (with $s_2 = 3$ and $t_2 = 4$) and
$\B(5, 4) = \{2$, $3\}$ (with $s_2 = s_3 = 2$ and $t_2 = t_3 = 4$).
These results show that
\begin{align*}
  c(5, 2)
  &=
  \Pi(5, 2) \,\,
  \exp( 4 B(5, 2) )
  \prod_{p \equiv 3 \bmod 5}
    \exp\Bigl( 4 f_{4, 3} \Bigl( \frac1p \Bigr) \Bigr)
  \exp\Bigl(
    4 \sum_{p \equiv 2 \bmod 5}
        \Bigl( f_{4, 1} \Bigl( \frac1p \Bigr) - \frac1p \Bigr)
    \Bigr) \\
  c(5, 3)
  &=
  \Pi(5, 3) \,\,
  \exp( 4 B(5, 3) )
  \prod_{p \equiv 2 \bmod 5}
    \exp\Bigl( 4 f_{4, 3} \Bigl( \frac1p \Bigr) \Bigr)
  \exp\Bigl(
    4 \sum_{p \equiv 3 \bmod 5}
        \Bigl( f_{4, 1} \Bigl( \frac1p \Bigr) - \frac1p \Bigr)
    \Bigr) \\
  c(5, 4)
  &=
  \Pi(5, 4) \,\,
  \exp( 4 B(5, 4) )
  \prod_{\substack{b \in \{2, 3 \} \\ p \equiv b \bmod 5}}
    \exp\Bigl( 4 f_{4, 2} \Bigl( \frac1p \Bigr) \Bigr)
  \exp\Bigl(
    4 \sum_{p \equiv 4 \bmod 5}
        \Bigl( f_{2, 1} \Bigl( \frac1p \Bigr) - \frac1p \Bigr)
    \Bigr),
\end{align*}
and the corresponding values for the constants $C(5, a)$ can be found
using \eqref{def-c}.
In the case of $C(5, 4)$ we can be slightly more explicit since
$f_{2, 1}(x) = 1 / 2 \log((1 + x) / (1 - x))$, so that
\[
  c(5, 4)
  =
  \Pi(5, 4) \,\,
  \exp( 4 B(5, 4) )
  \prod_{\substack{b \in \{2, 3 \} \\ p \equiv b \bmod 5}}
    \Bigl( \frac{1 + p^{-2}}{1 - p^{-2}} \Bigr)
  \prod_{p \equiv 4 \bmod 5}
    \Bigl( \frac{1 + p^{-1}}{1 - p^{-1}} \Bigr)^2 e^{- 4 / p}.
\]

\bibliographystyle{plain}

\begin{tabular}{l@{\hskip 14mm}l}
A.~Languasco               & A.~Zaccagnini \\
Universit\`a di Padova     & Universit\`a di Parma \\
Dipartimento di Matematica & Dipartimento di Matematica \\
Pura e Applicata           & Parco Area delle Scienze, 53/a \\
Via Trieste 63             & Campus Universitario \\
35121 Padova, Italy        & 43100 Parma, Italy \\
{\it e-mail}: languasco@math.unipd.it & {\it e-mail}:
alessandro.zaccagnini@unipr.it
\end{tabular}

\end{document}